\date{Rev. 5/VII/2011}
\title{The determinant bound for discrepancy is almost
tight}

\newcommand{\cmt}[1]{\ifhmode\newline\fi{\sf *** \ \ #1 \\}}

\author{
{\sc Ji\v{r}\'{\i} Matou\v{s}ek}\thanks{Partially supported
by the  ERC Advanced Grant No.~267165.
}
\\
   {\footnotesize Department of Applied Mathematics and}\\[-1.5mm]
   {\footnotesize Institute of Theoretical Computer Science (ITI)}\\[-1.5mm]
   {\footnotesize  Charles University, Malostransk\'{e} n\'{a}m. 25}\\[-1.5mm]
{\footnotesize  118~00~~Praha~1,
   Czech Republic}
}

\documentclass[11pt]{article}

\usepackage{a4wide}
\usepackage{latexsym}
\usepackage{amssymb}
\usepackage{epsfig}
\usepackage{url}

\newtheorem{theorem}{Theorem}

\newtheorem{lemma}[theorem]{Lemma}

\newcommand{\heading}[1]{\vspace{1ex}\par\noindent{\bf #1}}

\makeatletter
\newcommand{\ProofEndBox}{{\ifhmode\unskip\nobreak\hfil\penalty50 \else
          \leavevmode\fi\quad\vadjust{}\nobreak\hfill$\Box$
            \finalhyphendemerits=0 \par}}
\makeatother
\newcommand{\proofend}{\ProofEndBox\smallskip}

\newcommand{\R}{{\mathbb{R}}}

\newcommand{\sgn}{\mathop {\rm sgn}\nolimits}

\newcommand\FF{\mathcal{F}}

\newcommand\makevec[1]{{\bf #1}}
\def \uu {\makevec{u}}

\def \ww {\makevec{w}}
\def \aa {\makevec{a}}
\def \bb {\makevec{b}}

\def \ee {\makevec{e}}
\def \xx {\makevec{x}}
\def \yy {\makevec{y}}
\def \zz {\makevec{z}}
\def \cc {\makevec{c}}
\def \nula {\makevec{0}}
\def \bone {\makevec{1}}

\def\bgamma{\mbox{{\boldmath $\gamma$}}}
\def\bbeta{\mbox{{\boldmath $\beta$}}}
\def\bxi{\mbox{{\boldmath $\xi$}}}
\def\bomega{\mbox{{\boldmath $\omega$}}}

\newcommand\PSD{{\mathrm{PSD}}}
\newcommand\SYM{{\mathrm{SYM}}}

\newcommand\scalp[2]{\langle#1,#2\rangle}
\def\:{\colon}

\long\def\onefigure#1#2{
\begin{figure*}[tbp]
\begin{center}
#1
\end{center}
\caption{#2}
\end{figure*}
}

\newcommand{\labepsfig}[2]  
{\onefigure{\mbox{\epsfig{file=#1.eps}}}{\label{f:#1} #2} }

\newcommand{\labepsfigw}[3]  
{\onefigure{\mbox{\epsfig{file=#1.eps,width=#2}}}{\label{f:#1} #3} }

\newcommand\disc{{\mathop{\rm disc}\nolimits}}
\newcommand\detlb{{\mathop{\rm detlb}\nolimits}}
\newcommand\herdisc{{\mathop{\rm herdisc}\nolimits}}
\newcommand\vecdisc{{\mathop{\rm vecdisc}\nolimits}}
\newcommand\hervecdisc{{\mathop{\rm hervecdisc}\nolimits}}

\newcommand\notion[1]{\emph{#1}}

\begin{document}

\maketitle

\begin{abstract} In 1986 Lov\'asz, Spencer, and Vesztergombi proved
a lower bound for the hereditary discrepancy of a set system $\FF$
in terms of determinants of square submatrices of the incidence
matrix of $\FF$. As shown by an example of Hoffman, this bound
can differ from $\herdisc(\FF)$ by a multiplicative factor of
order almost $\log n$, where $n$ is the size of the ground
set of $\FF$. We prove that it never differs by more than
$O((\log n)^{3/2})$, assuming $|\FF|$ bounded by a polynomial in~$n$.
We also prove that if such an $\FF$ is the union of $t$ systems
$\FF_1,\ldots,\FF_t$, each of hereditary discrepancy at most $D$, 
then $\herdisc(\FF)\le
O(\sqrt t (\log n)^{3/2}D)$. For $t=2$, this 
almost answers a question of S\'os. 
The proof is based on a recent algorithmic result of Bansal,
which computes low-discrepancy colorings using semidefinite
programming.
\end{abstract}

\section{Introduction}

Let $V=[n]:=\{1,2,\ldots,n\}$ be a vertex set and $\FF=\{F_1,F_2,\ldots,F_m\}$
be a system of subsets of~$V$.
The \notion{discrepancy} of $\FF$ is
$\disc(\FF):=\min_\chi \disc(\FF,\chi)$,
where the minimum is over all colorings $\chi\:V\to\{-1,+1\}$,
and
$\disc(\FF,\chi):=\max_{i=1,2,\ldots,m}\bigl|\sum_{j\in F_i}\chi(j)\bigr|$.

The \notion{hereditary discrepancy} of $\FF$ is
$$
\herdisc(\FF):=\max_{J\subseteq V}\disc(\FF|_J).
$$
Here $\FF|_J$ denotes the \emph{restriction} of the set system $\FF$
to the ground set $J$, i.e., $\{F\cap J:F\in\FF\}$.

Bounding the discrepancy or the hereditary discrepancy of a particular
set system from below is usually challenging. One of the strongest
known tools is a result known
as the \notion{determinant lower bound}.
To formulate it we define, for a real matrix $A$, 
$$
\detlb(A) := \max_{k}\max_B |\det B|^{1/k},
$$
where the maximum is over all $k\times k$ submatrices $B$
of $A$. For a set system $\FF$, we put
$\detlb(\FF):=\detlb(A)$, where $A$ is the incidence matrix of~$\FF$.

\begin{theorem}[Lov\'asz, Spencer,
and Vesztergombi \cite{lsv-dssm-86}]\label{t:detlb}
For every (finite) set system~$\FF$ we
have\footnote{The bound in \cite{lsv-dssm-86}
is stated without the $\frac12$ factor. This has two causes:
first, their discrepancy is scaled by $\frac12$ compared to ours,
and second, in their argument,
at one step they seem to be multiplying by $2$ where, in my opinion,
one should divide by~2.} 
$\herdisc(\FF)\ge {\textstyle \frac 12} \detlb(\FF).
$
\end{theorem}

Lov\'asz et al.~\cite{lsv-dssm-86} conjectured that
the determinant lower bound is tight up to a constant factor,
i.e., $\herdisc(\FF)=O(\detlb(\FF))$ for all $\FF$. This was
refuted by an example of 
Hoffmann,\footnote{The vertex set in the example is the
set of edges of the complete $k$-ary tree $T$ of depth $k$
(so $n\approx k^k$). The set system $\FF$ is
a union $\FF_1\cup\FF_2$, where  $\FF_1$
consists of  the edge sets of all root-to-leaf paths,
and $\FF_2$ contains,
for each non-leaf vertex $v$, the set of the $k$ edges 
connecting $v$ to its successors. 
We have $\herdisc(\FF_1)\le 1$, $\herdisc(\FF_2)\le 1$,
$\detlb(\FF)=O(1)$, and 
$\disc(\FF)=k
\approx (\log n)/(\log\log n)$.
See, e.g., \cite{bs-dt-95} or \cite{m-gd}
for more details.} 
which shows that $\herdisc(\FF)/\detlb(\FF)$ can be
of order $(\log n)/(\log\log n)$. A construction
of P\'alv\"olgyi \cite{palvo-wedges}, also
presented in Section~\ref{s:dom} below, provides the slightly stronger
lower bound of $\Omega(\log n)$ for the same quantity.

Here we prove that $\herdisc(\FF)/\detlb(\FF)$
cannot be much larger than in these examples, at least
if $|\FF|$ is bounded by a polynomial function of~$n$.

\begin{theorem}\label{t:}
For every set system $\FF$
$$
\herdisc(\FF)\le \detlb(\FF)\cdot O\Bigl(\log(mn)\sqrt{\log n}\,\Bigr).
$$
\end{theorem}

Next, we consider the situation where a set system $\FF$
as above is a union of set systems $\FF_1,\FF_2,\ldots,\FF_t$,
and we are interested in bounding $\herdisc(\FF)$ in terms
of the hereditary discrepancies of the $\FF_i$.

For $t=2$, this problem was raised by S\'os 
(it is cited, e.g., in \cite{lsv-dssm-86}). She asked
whether $\herdisc(\FF_1\cup\FF_2)$ can be estimated in terms
of $\herdisc(\FF_1)$ and $\herdisc(\FF_2)$ for any two
set systems $\FF_1$ and $\FF_2$ (on the same vertex set).
Hoffman's example mentioned
above shows that $\herdisc(\FF_1\cup\FF_2)$ cannot be
bounded by a function of $\herdisc(\FF_1)$ and $\herdisc(\FF_2)$ alone.

The next theorem shows that a good bound is possible
if we also allow for a moderate dependence on $m$ and $n$.
Namely, $\herdisc(\FF_1\cup\FF_2)$ can exceed
$\max(\herdisc(\FF_1),\herdisc(\FF_2))$ at most by a factor
polylogarithmic in $n$ and $m$, not much more than in
Hoffman's example.\footnote{A connection of S\'os's question
to the tightness of the determinant lower bound
was observed in \cite{lsv-dssm-86}, although with
a different proof, which would yield a quantitatively weaker
result in our setting.}

The only previous result in this direction,
from \cite{kmv-dass}, shows that if $\FF_2$ consists of a single
set, then $\herdisc(\FF_1\cup\FF_2)=O(\herdisc(\FF_1)\log n)$.

\begin{theorem}\label{t:union}
Let $\FF$ be a system of $m$ sets on $n$ vertices,
and let $\FF=\FF_1\cup\FF_2\cup\cdots\cup\FF_t$.
Then
$$
\herdisc(\FF)\le \bigl(\max_{i=1,2,\ldots,t}\herdisc (\FF_i)\bigr)\cdot O\Bigl(\sqrt t\,\log(mn)\sqrt{\log n}\Bigr).
$$
\end{theorem}

There exist systems $\FF$ of $n$ sets on $n$ points with discrepancy
of order $\sqrt n$ (e.g., systems derived from Hadamard matrices
or random set systems; see, e.g., \cite{m-gd}). If we let
$\FF_i$ be the set system consisting of the $i$th set of
such an $n$, $i=1,2,\ldots,n$, 
then $\herdisc(\FF_i)=1$, while $\herdisc(\FF)=\Omega(\sqrt n)$.
In this sense, the bound in Theorem~\ref{t:union} is tight
up to a polylogarithmic factor, including the dependence on~$t$.

Theorem~\ref{t:union} is an immediate consequence of Theorems~\ref{t:detlb}
and~\ref{t:} and of the next linear-algebraic lemma.

\begin{lemma}\label{l:union}
Let $A_1,\ldots,A_t$ be real matrices, each with $n$ columns,
let $D:=\max_{i=1,2,\ldots,t} \detlb(A_i)$, and
let $A$ be a matrix in which each row is a row of some of the $A_i$.
Then
$$
\detlb(A)\le D\sqrt{et}.
$$
\end{lemma}

In case where some of the $\detlb(A_i)$ are much smaller
than the others, for example, one can obtain a somewhat
better upper bound by making finer estimates in the
calculation in the proof. However, a general formulation
of such a finer bound looks cumbersome, and it seems that
in such cases, a similar improvement can usually be achieved
by applying the lemma repeatedly in several stages
with various values of~$D$. 

The bound in Lemma~\ref{l:union} is generally tight up to a constant
factor;
this can be seen from an example similar to the one
mentioned below Theorem~\ref{t:union}. Namely,
let  $A$ be an $n\times n$ Hadamard matrix
(i.e., a matrix with pairwise orthogonal rows
and with  $\pm 1$ entries, which is well known
to exist for infinitely many values of~$n$), and let
$A_i$ be the single-row matrix made of the $i$th row of $A$,
$i=1,2,\ldots,n$. Then, obviously, $\detlb(A_i)=1$,
while $\detlb(A)\ge \det(A)^{1/n}=\sqrt n$.
Moreover, if we partition the rows of this $A$
into $t$ blocks $A_1,\ldots,A_t$ by $n/t$ rows each, then
the Hadamard bound, stating that the determinant
of a matrix is at most the product of the Euclidean
norms of the rows, implies $\detlb(A_i)\le \sqrt{n/t}$.
This shows the tightness of Lemma~\ref{l:union} 
for all~$t\le n$.

\heading{Vector discrepancy and Bansal's algorithm. } The proof
of Theorem~\ref{t:} is based on a recent breakthrough---an
algorithm of Bansal \cite{bansal-di}. The algorithm produces
a low-discrepancy coloring of a given set system,
using semidefinite programming and a clever randomized rounding
strategy. To state the consequence of Bansal's work
that we will use, we first introduce another notion of
discrepancy.

The \notion{vector discrepancy} of the set system $\FF$, denoted by
$\vecdisc(\FF)$, is the smallest $D\ge 0$ for which
there exist unit vectors $\uu_1,\ldots,\uu_n\in\R^n$
such that
$$
\biggl\|\sum_{j\in F_i} \uu_j\biggr\| \le  D,\ \ \ i=1,2,\ldots,m,
$$
where $\|\cdot\|$ is the Euclidean norm.
So, for vector discrepancy, one colors by unit vectors
instead of $\pm1$'s. Since a $\pm1$ coloring can also be regarded
as a vector coloring by the vectors $\ee_1=(1,0,\ldots,0)$ and $-\ee_1$,
we have $\vecdisc(\FF)\le\disc(\FF)$.

The  \notion{hereditary vector discrepancy}
$\hervecdisc(\FF)$ is the maximum vector discrepancy of
a restriction of $\FF$ to a subset $J\subseteq V$.

Bansal's algorithm yields the following.

\begin{theorem}[Bansal \cite{bansal-di}]\label{t:bans}
We have 
$$\herdisc(\FF)=O\bigl(\log(mn)\bigr)\cdot\hervecdisc(\FF)
$$ 
for all~$\FF$.
\end{theorem}

We conjecture that the claim of Theorem~\ref{t:bans}
actually holds with $\sqrt{\log (mn)}$ instead of $\log(mn)$.
If true, this would yield a similar improvement in
Theorem~\ref{t:} and get close to an asymptotically optimal
bound, at least assuming $m$ bounded by a polynomial in~$n$.

\section{A dual formulation of the vector discrepancy}


\begin{lemma}\label{l:dual}
We have $\vecdisc(\FF)\ge D$ if and only if there are
nonnegative reals $w_1,\ldots,w_m$ with $\sum_{i=1}^m w_i\le 1$ and
reals $z_1,\ldots,z_n$ with $\sum_{j=1}^n z_j\ge D^2$
such that for all $\xx\in\R^n$,
\begin{equation}\label{e:dui}
\sum_{i=1}^m w_i \biggl(\sum_{j\in F_i} x_j\biggr)^2 \ge
\sum_{j=1}^n z_j x_j^2.
\end{equation}
\end{lemma}

\heading{Proof. } 
We will use the \emph{duality of semidefinite programming}.
Dualizing a semidefinite program is a routine procedure,
but unfortunately, I am not aware of an explicit recipe for the general
case in the literature. Rather than converting the relevant semidefinite
program to a standard form, it seems more convenient to
use a duality theorem for conic programing
from Duffin~\cite{Duf}, which we now introduce.

Let $V,W$ be real vector spaces
(for our purposes we may assume that they are finite-dimensional), 
each with a scalar product, which we denote by $\scalp\cdot\cdot$
in both cases. Let $K\subseteq V$ and
$L\subseteq W$ be closed convex cones,\footnote{A \emph{convex
cone} is a convex set $K$ such that $\xx\in K$ implies
$\lambda\xx\in K$ for all $\lambda\ge 0$. The \emph{dual cone}
of $K$ is $K^*=\{\yy:\scalp\xx\yy\ge 0\mbox{ for all }\xx\in K\}$.
A simple property we will often use is $(K\oplus L)^*=
K^*\oplus L^*$, where $\oplus$ denotes direct sum;
we assume $K\subseteq V$, $L\subseteq W$, where $V,W$ are disjoint
vector spaces, and $K\oplus L\subseteq V\oplus W$.}
let $\bb\in W$ and $\cc\in V$ be vectors, and 
let $F\:V\to W$ be a linear map. We consider the \emph{primal cone
program} (P)
$$
\mbox{maximize } \scalp\cc\bxi\mbox{ subject to } \bb-F(\bxi)\in L,\bxi\in K
$$
and the \emph{dual cone program} (D)
$$
\mbox{minimize } \scalp\bb\bomega\mbox{ subject to } F^T(\bomega)-\cc\in K^*,
\bomega\in L^*.
$$
Here $F^T\:W\to V$ denotes the \emph{adjoint} of $F$; if we fix orthonormal
bases in
$V$ and $W$, then the matrix of $F^T$ is the transpose of the matrix
of $F$. The duality theorem asserts that if the maximum in (P)
is a finite number $\gamma$ and if the set of feasible solutions
of (P) has an interior point, then (D) is feasible as well
and its minimum equals~$\gamma$.  (Since (P) and (D) are dual to one
another, one can also interchange their role in the theorem.)

In our case, we start with a vector program defining $\vecdisc(\FF)$,
namely,
$$
\mbox{minimize }t\mbox{ subject to }
\|\textstyle\sum_{j\in F_i} \uu_j\|^2\le t, i=1,2,\ldots,m,
\|\uu_1\|=\cdots=\|\uu_n\|=1,
$$
and we convert it to an equivalent semidefinite program
in a standard way. We introduce a variable $Q$, which
is an $n\times n$ matrix with $q_{ij}=\uu_i^T\uu_j$ (the Gram
matrix of the $\uu_j$); as is well known, $Q$ is of this form
for some vectors $\uu_1,\ldots,\uu_n$ exactly if $Q\in\PSD_n$,
where $\PSD_n$ denotes the cone of $n\times n$
positive semidefinite matrices. 

For the constraint $\|\sum_{j\in F_i} \uu_j\|^2\le t$,
we expand the left-hand side to $\sum_{j,k\in F_i}\uu_j^T\uu_k$,
and then it translates to
\begin{equation}\label{e:co1}
t-(\aa_i\aa_i^T)\bullet Q\ge 0,
\end{equation}
where $\aa_i$ is the $i$th row of the incidence matrix $A$ of $\FF$,
regarded as an $n\times 1$ matrix (so $\aa_i\aa_i^T$ is an $n\times n$
matrix), and $\bullet$ denotes the scalar product of matrices
(given by
$X\bullet Y=\sum_{i,j} x_{ij}y_{ij}$). The constraint $\|\uu_j\|=1$
then reads
\begin{equation}\label{e:co2}
E_j\bullet Q=1,
\end{equation}
where $E_j$ is the matrix with $1$ at position $(j,j)$ and $0$s elsewhere.

This semidefinite program can be regarded as the cone program (P)
with $V=\SYM_n\oplus \R$ (where $\SYM_n$ denotes the vector
space of all symmetric $n\times n$ matrices with the $\bullet$
scalar product), $W=\R^m\oplus\R^n$, $K=\PSD_n\oplus[0,\infty)$, and
$L=[0,\infty)^m\oplus \{\nula_n\}$
(where $\nula_n$ denotes the $n$-component vector of $0$s). 
We write the unknown
$\bxi\in V$ in the form $(Q,t)$; then we have $F(Q,t)=
(\aa_1\aa_1^T\bullet Q-t,\ldots,\aa_m\aa_m^T\bullet Q-t,
E_1\bullet Q,\ldots,E_n\bullet Q)$. Finally, $\cc=(0_{n\times n},-1)$
and $\bb=(\nula_m,\bone_n)$.

As for the dual cone program, it is well known that $\PSD_n^*=\PSD_n$,
so $K^*=K$,
and it is easily checked that $L^*=[0,\infty)^m\oplus \R^n$. Let us write
the variable $\bomega$ in (D) in the form $(\ww,-\zz)$,
$\ww\in\R^m$, $\zz\in\R^n$. Then (D) becomes
\begin{equation}\label{thisp}
\begin{array}{rl}
\mbox{ maximize } \textstyle \sum_{j=1}^n z_j
\mbox{ subject to} & \textstyle\sum_{i=1}^m w_i \aa_i\aa_i^T -\sum_{j=1}^n z_jE_j\in\PSD_n,\\
&\ \ \textstyle1-\sum_{i=1}^m w_j\ge 0, w_1,\ldots,w_m\ge 0.
\end{array}
\end{equation}
Let $M$ stand for the matrix 
$\sum_{i=1}^m w_i \aa_i\aa_i^T -\sum_{j=1}^n z_jE_j$.
Then the condition of positive semidefiniteness of $M$
means $\xx^T M\xx\ge 0$ for all $\xx\in\R^n$,
and it is easy to verify that this
can be rewritten as the inequality (\ref{e:dui})
in the lemma.

It remains to verify that the duality theorem can be applied
to these (P) and (D). We will check that (D) is bounded and
has a feasible interior point. To verify boundedness, which
means that 
$\sum_{j=1}^nz_j$ cannot be arbitrarily large,
we use (\ref{e:dui}) with $\xx=\bone_n$: then the left-hand
side is bounded, and so $\sum_{j=1}^nz_j$ is bounded as well.

For an interior feasible point, we can take, e.g.,
$w_1=\cdots=w_m=\frac1{2m}$, $z_1=\cdots=z_n=-1$.
Then (\ref{e:dui}) obviously holds for all $\xx$.

Thus, the duality theorem applies and shows that 
if $\vecdisc(\FF)\ge D$,
then the maximum in (\ref{thisp}) is at least $D^2$.
This yields the existence of the desired $w_i$ and $z_j$, and
 the lemma is proved.
\proofend


\section{Proof of Theorem~\ref{t:}}

We begin with a simple and probably standard lemma.

\begin{lemma}\label{l:almostconst}
 Let $\yy\in \R^n$ be a vector. Then
there exists a subset $K\subseteq [n]$ of indices
such that $\yy$ is ``almost constant'' on $K$,
in the sense that, for some $t>0$, we have
 $t< |y_j|\le 2t$ for all $j\in K$, and 
$$
\|\yy[K]\| \ge \Omega\left(\frac{\|\yy\|}{\sqrt{\log n}}\right),
$$
where
$\yy[K]$ denotes the $|K|$-component vector
$(y_j:j\in K)$.
\end{lemma}

\heading{Proof. }
Let $y_{\rm max} :=
\max_j|y_j|$, and for $i=0,1,2,\ldots$, let
$K_i :=\{j: |y_j|\in (2^{-i-1}y_{\rm max}, 
2^{-i}y_{\rm max}]\}$. The contribution to $\|\yy\|$ of the components
of $\yy$ with indices in $K_i$ for $i\ge 2\log n$, say,
is negligible, and so there exists some $i_0$ for which
$\sum_{j\in K_{i_0}}y_j^2 =\Omega(\|\yy\|^2/\log n)$.
Then $K:=K_{i_0}$ will do.
\proofend

Theorem~\ref{t:} will follow from Bansal's result
(Theorem~\ref{t:bans}) and the next lemma.

\begin{lemma}\label{l:main}
Let $\FF=\{F_1,\ldots,F_m\}$ be a set system  on $[n]$
with $\vecdisc(\FF)=D$. Then 
$\detlb(\FF)= \Omega(D/\sqrt{\log n}\,)$.
\end{lemma}

\heading{Proof. } We begin with the dual formulation
of vector discrepancy from Lemma~\ref{l:dual}.
For more convenient notation, we will write
the nonnegative weight $w_i$ as $\beta_i^2$.
Moreover, we let $J\subseteq [n]$ consist of the indices
$j$ with $z_j>0$, and we will use the inequality (\ref{e:dui})
in Lemma~\ref{l:dual} only for vectors $\xx$ that are zero
outside $J$. Writing $z_j=\gamma_j^2$ for $j\in J$,
we arrive at the inequality
\begin{equation}\label{e:dual1}
\sum_{i=1}^m \beta_i^2\biggl(\sum_{j\in F_i\cap J} x_j\biggr)^2
\ge \sum_{j\in J}^n \gamma_j^2 x_j^2
\end{equation}
for all $\xx\in\R^J$, where $\|\bbeta\|\le 1$ and
$\|\bgamma\|\ge D$.


Next,
using Lemma~\ref{l:almostconst} for $\yy=\bgamma$, we select
$K\subseteq J$ with $\|\bgamma[K]\|= \Omega(D/\sqrt{\log n})$
and with $\bgamma[K]$ almost constant (within a factor
of 2). Setting $k:=|K|$ and $\tilde D:=\frac 12\|\bgamma[K]\|$,
the quadratic average
of $\gamma_j$ over $j\in K$ equals $2\tilde D/\sqrt k$,
and so $\gamma_j\ge \tilde D/\sqrt k$ 
for all $j\in K$. Therefore, restricting (\ref{e:dual1})
to vectors $\xx$ with $x_j=0$ for $j\not\in K$,
we obtain
\begin{equation}\label{e:2}
\sum_{i=1}^m \beta_i^2\biggl(\sum_{j\in F_i\cap K} x_j\biggr)^2
\ge \frac{\tilde D^2}{k}\sum_{j\in K} x_j^2
\end{equation}
for all $\xx\in\R^K$. 

Let $C:=A[*,K]$ be the $m\times k$ incidence matrix of the system $\FF|_K$
(consisting of the columns of $A$ whose indices belong to $K$),
and let $\check C$ be the $m\times k$ matrix obtained from $C$
by multiplying the $i$th row by $\beta_i$.
Then (\ref{e:2}) can be rewritten as
$$
\xx^T\check C^T\check C\xx=\|\check C\xx\|^2\ge \frac{\tilde D^2}{k}\|\xx\|
$$
for all $\xx\in\R^k$. 

This, by the usual variational characterization
of eigenvalues, tells us that the smallest eigenvalue of 
the $k\times k$ matrix $\check C^T\check C$ satisfies
$\lambda_{\rm min}(\check C^T \check C)\ge \tilde D^2/k$.
Then, since the determinant is the product of eigenvalues,
we also have $\det(\check C^T \check C)\ge (\tilde D^2/k)^k$.
The Binet--Cauchy formula then asserts that
\begin{equation}\label{e:3}
\det(\check C^T \check C)=\sum_I \det(\check C[I,*])^2,
\end{equation}
where the summation is over all $k$-element subsets
$I\subseteq [m]$ and $\check C[I,*]$ consists of
the rows of $\check C$ whose indices lie in $I$.

We have $\det(\check C[I,*])=\det(C[I,*])\prod_{i\in I}\beta_i$.
Setting $M:=\max_I |\det(C[I,*])|$, we can rewrite
the right-hand side of (\ref{e:3}) and estimate it as follows:
\begin{eqnarray*}
\sum_I \det(\check C[I,*])^2&=&
\sum_I \det(C[I,*])^2\prod_{i\in I}\beta_i^2\\
&\le& M^2\sum_{I} \prod_{i\in I}\beta_i^2\\
&\le& M^2\frac{\left(\sum_{i=1}^m \beta_i^2\right)^k}{k!} \le
\frac{M^2}{k!},
\end{eqnarray*} 
where the penultimate inequality follows because
every term $\prod_{i\in I}\beta_i^2$ occurs $k!$ times in the
multinomial expansion of $(\beta_1^2+\cdots+\beta_m^2)^k$.

Letting $B:=C[I,*]$ for an $I$ maximizing $|\det C[I,*]|$,
we thus have
$\det(B)^2 \ge  k! \det(\check C^T \check C)\ge k!(\tilde D^2/k)^k
\ge (k/e)^k (\tilde D^2/k)^k=\Omega(\tilde D)^{2k}=
\Omega(D/\sqrt{\log n})^{2k}$. So the $k\times k$ matrix $B$
witnesses $\detlb(\FF)= \Omega(D/\sqrt{\log n}\,)$, and the lemma
is proved.
\proofend

\heading{Proof of Theorem~\ref{t:}. }
By Theorem~\ref{t:bans}, there is a subset $J\subseteq [n]$
with $\vecdisc(\FF|_J)=\Omega(\herdisc(\FF)/\log(mn))$.
Theorem~\ref{t:} follows by applying Lemma~\ref{l:main}
to $\FF|_J$.
\proofend

\section{Proof of Lemma~\ref{l:union}}



Let us consider a $k\times k$ submatrix $B$ of the matrix
$A$, and let $I_1,\ldots,I_t$ be index sets partitioning
$[k]$ such that $B[I_\ell,*]$ is a submatrix of $A_\ell$,
$\ell=1,2,\ldots,t$ (we also admit $I_\ell=\emptyset$).
Let $k_\ell:=|I_\ell|$.

We now apply the Gram--Schmidt orthogonalization to the row
vectors of each $B[I_\ell,*]$, separately for each $\ell$. 
We obtain a new matrix $\tilde B$ 
such that the rows of each $\tilde B[I_\ell,*]$
are orthogonal vectors.  Let $\tilde\bb_i$ be the
$i$th row of $\tilde B$.

Since the orthogonalization proceeds
by elementary row operations, which do not change the determinant,
we have $\det\tilde B=\det B$. We estimate $\det B$
using the Hadamard bound
\begin{equation}\label{e:1}
|\det B|=|\det\tilde B|\le\prod_{i=1}^k \|\tilde\bb_i\|=\prod_{\ell=1}^t
\prod_{i\in I_\ell}\|\tilde\bb_i\|.
\end{equation}

Let us fix $\ell$ for a moment.
In order to estimate $\prod_{i\in I_\ell}\|\tilde\bb_i\|$,
we consider the quantity 
$\det(\tilde B[I_\ell,*] \tilde B[I_\ell,*]^T)$.
On the one hand, since the $\tilde\bb_i$ are mutually
orthogonal for $i\in I_\ell$, the matrix $\tilde B[I_\ell,*] \tilde B[I_\ell,*]^T$ is diagonal with the entries $\|\tilde\bb_i\|^2$, $i\in I_\ell$,
on the diagonal, and so 
\begin{equation}\label{e:a2}
\det (\tilde B[I_\ell,*]\tilde B[I_\ell,*]^T)=
\prod_{i\in I_\ell} \|\tilde\bb_i\|^2.
\end{equation}

On the other hand, by the Binet--Cauchy formula, we have
\begin{equation}\label{e:a3}
\det(\tilde B[I_\ell,*] \tilde B[I_\ell,*]^T)=
\sum_{J\subseteq[k],|J|=k_\ell}\det(\tilde B[I_\ell,J])^2=
\sum_{J\subseteq[k],|J|=k_\ell}\det( B[I_\ell,J])^2,
\end{equation}
the last equality holding again because $\tilde B[I_\ell,J]$
is obtained from $B[I_\ell,J]$ by elementary row operations
which preserve the determinant. 

Putting (\ref{e:1}), (\ref{e:a2}), and (\ref{e:a3}) together, we arrive at
$$
\det(B)^2 \le 
\prod_{\ell=1}^t \sum_{J\subseteq[k],|J|=k_\ell}\det( B[I_\ell,J])^2
\le \prod_{\ell=1}^t {k\choose k_\ell}\detlb(A_\ell)^{2k_\ell}
\le D^{2k} \prod_{\ell=1}^t {k\choose k_\ell},
$$
since each $B[I_\ell,J]$ is a $k_\ell\times k_\ell$ submatrix of~$A_\ell$.

Then we estimate, using the concavity of the function
$x\mapsto x\ln (1/x)$ and Jensen's inequality,
\begin{eqnarray*}
\prod_{\ell=1}^t {k\choose k_\ell}&\le& \prod_{\ell=1}^t
\left(\frac{ek}{k_\ell}\right)^{k_\ell}=
e^{\sum_{\ell=1}^t k_\ell\ln(ek/k_\ell)}\\
&\le& e^{t(k/t)\ln(et)}=(et)^k.
\end{eqnarray*}
Thus $\det(B)^{1/k}\le D\sqrt{et}$ as claimed.
\proofend

\heading{Remark. } The case $t=2$ has a somewhat simpler
proof using  the \emph{Laplace expansion} of $\det B$, which
asserts that
$$
\det B = \sum_J \sgn(I,J) \det(B[I,J])\det(B[\overline I,
\overline J]),
$$
where the sum is over all $|I|$-element subsets $J\subseteq [k]$,
$\overline I=[k]\setminus I$,
and $\sgn(I,J)\in\{\pm 1\}$ is
a sign depending on  $I$ and $J$ in a way that is of no concern
for us (see, e.g., \cite[Theorem~4.3]{BhattaLA}).

In our case, we again let $I_\ell$ be the set of indices 
of the rows of $B$ that come from $A_\ell$, $\ell=1,2$,
and $k_\ell=|I_\ell|$.
We use the Laplace expansion of $\det B$ for $I=I_1$.
For every $J$ we have $|\det(B[I_1,J])|\le D^{k_1}$
and $|\det(B[I_2,\overline J])|\le D^{k_2}$ by the assumption, and so
$\det B \le {k\choose k_1} D^k\le 2^kD^k$.
Hence $\det(B)^{1/k}\le 2D$ as needed.

A bound on $\detlb(A)$ for larger $t$ can also be
obtained by iterating this argument, but this method
apparently leads only to $\detlb(A)=O(tD)$.

\section{P\'alv\"olgyi's example}\label{s:dom}

As was pointed out by P\'alv\"olgyi (private communication, 2011),
his geometric construction in \cite{palvo-wedges}
actually provides a slight quantitative improvement
over Hoffman's example. Translated to the setting
of set systems, the construction yields, for every
$k\ge 1$, two systems $\FF_1,\FF_2$ of $k$-element
subsets of $[n]$, with $n={2k\choose k}-1< 4^k$,
such that 
\begin{enumerate}
\item[(i)] $\herdisc(\FF_1),\herdisc(\FF_2)\le 1$, and
\item[(ii)] under every two-coloring of $[n]$,
$\FF_1\cup\FF_2$ contains a monochromatic set (and consequently,
$\disc(\FF_1\cup\FF_2)=k$).
\end{enumerate}

Since the construction in \cite{palvo-wedges} is presented geometrically,
and property (i) is not entirely obvious, we provide a short self-contained
exposition. 

The construction is inductive, and it requires two parameters,
$k$ and $\ell$. The inductive hypothesis is the following:
\begin{quote}
For every $k,\ell\ge 1$ and a ground set $V$ of $n={k+\ell\choose k}-1$
elements,
there exist set systems
$\FF_1=\FF_1(V,k,\ell)$, $\FF_2=\FF_2(V,k,\ell)$ on $V$ 
such that $\FF_1$ consists of $k$-tuples, $\FF_2$ consists
of $\ell$-tuples, $\herdisc(\FF_1),\herdisc(\FF_2)\le 1$, and
for every red-blue coloring of $V$, there exists a set 
$F\in\FF_1$ that is completely red or a set $F\in\FF_2$ that
is completely blue. 
\end{quote}

For $k=1$ and $\ell$ arbitrary, we have $n=\ell$, and we 
take $\FF_1$ consisting of the $\ell$ singleton subsets of $V$,
while $\FF_2=\{V\}$. 
The other base case $\ell=1$ and $k$ arbitrary is symmetric.

In the inductive step, we take the ground set $V$ and decompose
it into three disjoint subsets: $V'$ with ${k+\ell-1\choose k-1}-1$
elements, $V''$ with ${k+\ell-1\choose k}-1$ elements, and
a singleton set $\{p\}$. We inductively construct
 $\FF_1':=\FF_1(V',k-1,\ell)$
and $\FF_2':=\FF_2(V',k-1,\ell)$, as well as
 $\FF_1'':=\FF_1(V'',k,\ell-1)$
and $\FF_2'':=\FF_2(V'',k,\ell-1)$.
Then we set $\FF_1:=\{F\cup\{p\}:F\in\FF_1'\}\cup\FF_1''$,
$\FF_2:=\FF_2'\cup\{F\cup\{p\}:F\in\FF_2''\}$, and check the
required properties by a straightforward induction.

We begin with verifying that, under every red-blue
coloring of $V$, $\FF_1$ has a completely red set or $\FF_2$
has a completely blue set. If $p$ is red,
then we restrict the coloring on $V'$ and inductively find
a red $(k-1)$-tuple in $\FF_1'$, which together with $p$
gives a red $k$-tuple in $\FF_1$, or a blue $\ell$-tuple
in $\FF_2'$, which is also a blue $\ell$-tuple in $\FF_2$.
The case of $p$ blue is analogous.

It remains to verify that $\herdisc(\FF_1),\herdisc(\FF_2)\le 1$.
By symmetry of the construction, it suffices to check
$\herdisc(\FF_1)$.
We inductively prove a slightly stronger statement:
for every $W\subseteq V$, there is a $\pm 1$ coloring $\chi$ of $W$
such that $\chi(F\cap W)\in\{0,1\}$ for every $F\in\FF_1$.
The base cases with $k=1$ or $\ell=1$ are obvious,
and in the inductive step, we inductively color $W'=V'\cap W$
and $W''=V''\cap W$, we flip the colors on $W'$, and we color $p$ 
with $+1$.

\subsection*{Acknowledgment}

I would like to thank Nikhil Bansal for enlightening e-mail
discussions concerning his algorithm and for help with 
dualizing the semidefinite program, and
D\"om\"ot\"or P\'alv\"olgyi for explaining me
how his construction in  \cite{palvo-wedges}
improves on Hoffmann's example.


\begin{thebibliography}{KMV05}

\bibitem[Ban10]{bansal-di}
N.~Bansal.
\newblock Constructive algorithms for discrepancy minimization.
\newblock \url{http://arxiv.org/abs/1002.2259}, also in \emph{FOCS'10: Proc.
  51st IEEE Symposium on Foundations of Computer Science}, pages 3--10, 2010.

\bibitem[BJN83]{BhattaLA}
P{.\,}B. Bhattacharya, S{.\,}K. Jain, and S{.\,}R. Nagpaul.
\newblock {\em {First course in linear algebra}}.
\newblock Wiley Eastern Limited, New Delhi etc., 1983.

\bibitem[BS95]{bs-dt-95}
J.~Beck and V.~{S\'os}.
\newblock Discrepancy theory.
\newblock In {\em Handbook of Combinatorics}, pages 1405--1446. North-Holland,
  Amsterdam, 1995.

\bibitem[Duf56]{Duf}
R{.\,}J. Duffin.
\newblock Infinite programs.
\newblock In H{.\,}W. Kuhn and A{.\,}W. Tucker, editors, {\em Linear
  Inequalities and Related Systems}, volume~38 of {\em Annals of Mathematical
  Studies}, pages 157--170. 1956.

\bibitem[KMV05]{kmv-dass}
J.-H. Kim, J.~Matou\v{s}ek, and V{.\,H.} Vu.
\newblock Discrepancy after adding a single set.
\newblock {\em Combinatorica}, 25:499--501, 2005.

\bibitem[LSV86]{lsv-dssm-86}
L.~Lov\'asz, J.~Spencer, and K.~Vesztergombi.
\newblock Discrepancy of set-systems and matrices.
\newblock {\em European J. Combin.}, 7:151--160, 1986.

\bibitem[Mat10]{m-gd}
J.~Matou\v{s}ek.
\newblock {\em Geometric Discrepancy (An Illustrated Guide), 2nd printing}.
\newblock Springer-Verlag, Berlin, 2010.

\bibitem[P{\'a}l10]{palvo-wedges}
D.~P{\'a}lv{\"o}lgyi.
\newblock {Indecomposable coverings with concave polygons}.
\newblock {\em Discrete Comput. Geom.}, 44(3):577--588, 2010.

\end{thebibliography}

\end{document}